\input amstex
\documentstyle{amsppt}
\magnification=\magstep1 \TagsOnLeft \NoBlackBoxes
\hsize=6.3 true in
\vsize=8.4 true in
\pageno=1
\TagsOnLeft

\def\rn{\Bbb R^n}
\def\bn{\Bbb R^{1+n}_+}

\topmatter
\title
Affine Variant of Fractional Sobolev Space with Application to Navier-Stokes System
\endtitle
\author
Jie Xiao  
\endauthor
\affil
\endaffil
\leftheadtext{Jie Xiao}
\rightheadtext{Affine Variant of Fractional Sobolev Space}

\thanks
\flushpar Research supported by NSERC (Canada) and Dean of Science Startup Fund (MUN, Canada). \flushpar 2000
\flushpar {\it Key words and phrases}: $W^{1,n}(\rn)$, $Q_\alpha(\rn)$, $BMO(\rn)$, Quadratic Morrey 
space, Sharp Sobolev imbedding, Incompressible Navier-Stokes system
\endthanks
\abstract
It is proved that for $\alpha\in (0,1)$, $Q_\alpha(\rn)$, not only as an intermediate space of $W^{1,n}(\rn)$ and $BMO(\rn)$ but also as an affine variant of Sobolev space $\dot{L}^{2}_\alpha(\rn)$ which is sharply imbedded in $L^{\frac{2n}{n-2\alpha}}(\rn)$, is isomorphic to a quadratic Morrey space under fractional differentiation. At the same time, the dot product $\nabla\cdot\big(Q_\alpha(\rn)\big)^n$ is applied to derive the well-posedness of the scaling invariant mild solutions of the incompressible Navier-Stokes system in $\bn=(0,\infty)\times\rn$.
\endabstract
\endtopmatter
\document

\heading
1. Introduction and Summary
\endheading

We begin by the square form of John-Nirenberg's BMO space (cf. [JN]) which plays an important role in harmonic analysis and applications to partial differential equations. For a locally integrable complex-valued function $f$ defined on the Euclidean space $\rn$, $n\ge 2$, with respect to the Lebesgue measure $dx$, we say that $f$ is of BMO class, denoted $f\in BMO=BMO(\rn)$, provided 
$$
\|f\|_{BMO}=\left(\sup_{I}{\big(\ell(I)\big)^{-n}}\int_{I}\big|f(x)-f_I\big|^2dx\right)^\frac12<\infty.
$$  
Here and elsewhere $\sup_I$ means that the supremum ranges over all cubes $I\subset\rn$ with edges parallel to the coordinate axes; $\ell(I)$ is the sidelength of $I$; and $f_I=\big(\ell(I)\big)^{-n}\int_I f(x)dx$ stands for the mean value of $f$ over $I$.

On the basis of the semi-norm $\|\cdot\|_{BMO}$, a large scale of function spaces has been introduced in [EsJPX], as defined below.

\proclaim{Definition 1.1} For $\alpha\in (-\infty,\infty)$, let $Q_\alpha$ be the space of all measurable complex-valued functions $f$ on $\rn$ obeying 
$$
\|f\|_{Q_\alpha}=\sup_{I}\left(\big(\ell(I)\big)^{2\alpha-n}\int_I\int_I\frac{|f(x)-f(y)|^2}{|x-y|^{n+2\alpha}}dxdy\right)^\frac12<\infty.
$$
\endproclaim

This $Q_\alpha$ is a natural extension of $BMO$ according to the following result (proved in [EsJPX] and [X]):
$$
Q_\alpha=\cases BMO,\quad &\alpha\in (-\infty,0),\\
\hbox{New\ space\ between}\ W^{1,n}\ \ \hbox{and}\ BMO,\quad &\alpha\in [0,1),\\
\Bbb C, \quad & \alpha\in [1,\infty).
\endcases
$$
Here $W^{1,n}=W^{1,n}(\rn)$ is the affine energy space of all $C^1$ functions $f$ on $\rn$ with
$$
\|f\|_{W^{1,n}}=\left(\int_{\rn}|\nabla f(x)|^n dx\right)^\frac{1}{n}<\infty.
$$
More importantly, $Q_\alpha$, $\alpha\in (0,1)$, can be regarded as the affinely invariant counterpart of the homogeneous Sobolev space $\dot{L}^2_\alpha=\dot{L}^2_\alpha(\rn)$ which consists of all complex-valued functions $f$ on $\rn$ with the $\alpha$-energy 
$$
\|f\|_{\dot{L}^2_\alpha}=\left(\int_{\rn}\int_{\rn}\frac{|f(x)-f(y)|^2}{|x-y|^{n+2\alpha}}dxdy\right)^\frac12<\infty.
$$
The reason for saying this is at least that $\|\cdot\|_{Q_\alpha}$ and $\|\cdot\|_{\dot{L}^2_\alpha}$ enjoy the following property:
$$
\|f\circ\phi\|_{Q_\alpha}=\|f\|_{Q_\alpha}\ \ \hbox{and}\ \ \|f\circ\phi\|_{\dot{L}^2_\alpha}=\lambda^{\alpha-\frac n2}\|f\|_{\dot{L}^2_\alpha}
$$
for any affine map $x\mapsto \phi(x)=\lambda x+ x_0$; $\lambda>0, x_0\in\rn$.

In the five year period since the paper [EsJPX] appeared, it has been found that $Q_\alpha$ is a useful and interesting concept; see also [AdL], [AnC], [PoS], [DX1], [DX2], [La], [PeY], [CuY] and [En]. This means that the study of this new space has not yet ended up -- in fact, there are many unexplored problems related to $Q_\alpha$. In this paper, although not attacking one of those open problems in Section 8 in [EsJPX], we go well beyond the previous results by studying the relation between this space and the quadratic Morrey space, but also giving an application of the induced facts to the incompressible Navier-Stokes system. 

To deal with the former, it is necessary to consider the following variant of [DX1, Theorem 3.3] that expands Fefferman-Stein's basic result for $BMO$ in [FS]: Given a $C^\infty$ function $\psi$ on $\rn$ with  
$$
\psi\in L^1,\quad |\psi(x)|\lesssim (1+|x|)^{-(n+1)},\quad\int_{\rn}\psi(x)dx=0\quad\hbox{and}\quad\psi_t(x)=t^{-n}{\psi(\frac{x}{t})}.
$$
Then for a measurable complex-valued function $f$ on $\rn$,
$$
f\in Q_\alpha\Longleftrightarrow\sup_{x\in\rn, r\in (0,\infty)}r^{2\alpha-n}\int_0^r\int_{|y-x|<r}{|f\ast\psi_t(y)|^2}{t^{-(1+2\alpha)}}dydt<\infty.\tag 1.1
$$
Here and henceforth, $L^p=L^p(\rn)$ represents the Lebesgue space equipped with $p$-norm $\|\cdot\|_{L^p}$; $\ast$ stands for the convolution operating on the space variable; and $U\lesssim V$ means that there exists a constant $c>0$ such that $U\le cV$. 

Two particular choices of $\psi$ in (1.1) yield two characterizations of $Q_\alpha$ involving the Poisson and heat semi-groups. As for this aspect, denote by $e^{-t\sqrt{-\Delta}}(\cdot,\cdot)$ and $e^{t\Delta}(\cdot,\cdot)$ are the Poisson and heat kernels respectively; that is, 
$$
e^{-t\sqrt{-\Delta}}(x,y)={\Gamma\big(\frac{n+1}{2}\big)}{\pi^{-\frac{n+1}{2}}}t{(|x-y|^2+t^2)^{-\frac{n+1}{2}}}
$$ 
and
$$
e^{t\Delta}(x,y)
=(4\pi t)^{-\frac{n}{2}}{\exp\Big(-\frac{|x-y|^2}{4t}\Big)}.
$$
Of course, for $\beta\in\Bbb R$, the notation $(-\Delta)^{\frac{\beta}{2}}$ is the $\beta/2$-th power of the Laplacian operator 
$$
-\Delta=-\sum_{j=1}^n\partial_j^2=-\sum_{j=1}^n\frac{\partial^2}{\partial x_j^2}
$$ 
determined by the Fourier transform $\widehat{(\cdot)}$
$$
\widehat{(-\Delta)^{\frac{\beta}{2}}f}(x)=(2\pi|x|)^\beta\hat{f}(x).
$$
On the one hand, if
$$
\psi_0(x)=\frac{1+|x|^2-(n+1)\Gamma\big(\frac{n+1}{2}\big){\pi^{-\frac{n+1}{2}}}}{(1+|x|^2)^{\frac{n+3}{2}}},
$$
then 
$$
(\psi_0)_t(x)=t{\partial_t}e^{-t\sqrt{-\Delta}}(x,0)
$$
and hence for $\alpha\in (0,1)$ and a measurable complex-valued function $f$ on $\rn$, 
$$
f\in Q_\alpha\Longleftrightarrow\sup_{x\in\rn, r\in (0,\infty)}r^{2\alpha-n}\int_0^r\int_{|y-x|<r}{\Big|{\partial_t}e^{-t\sqrt{-\Delta}}f(y)\Big|^2}{t^{1-2\alpha}}dydt<\infty.\tag 1.2
$$ 
On the other hand, if 
$$
\psi_j(x)=-(4\pi)^{-\frac n2}\Big(\frac{x_j}{2}\Big)\exp\Big(-\frac{|x|^2}{4}\Big)\quad\hbox{for}\quad j=1,...,n,
$$
then 
$$
(\psi_j)_t(x)=t{\partial_j}e^{t^2\Delta}(x,0)
$$
and so, for $\alpha\in (0,1)$ and a measurable complex-valued function $f$ on $\rn$,
$$
f\in Q_\alpha\Longleftrightarrow\sup_{x\in\rn, r\in (0,\infty)}r^{2\alpha-n}\int_0^r\int_{|y-x|<r}|\nabla e^{t^2\Delta}f(y)|^2{t^{1-2\alpha}}dydt<\infty.\tag 1.3
$$ 

With the help of the above-mentioned facts, we can establish the following result.

\proclaim{Theorem 1.1} Let $\alpha\in (0,1)$. Then

\flushpar {\rm (i)} $Q_\alpha=(-\Delta)^{-\frac{\alpha}{2}}\Cal L_{2,n-2\alpha}\hookrightarrow BMO$ is proper with
$$
\sup_{\|f\|_{Q_\alpha}>0}\frac{\|f\|_{BMO}}{\|f\|_{Q_\alpha}}\le\sqrt{\frac{n^{\frac{n+2\alpha}{2}}}{2}},
$$
where a measurable complex-valued function $f$ on $\rn$ belongs to $\Cal L_{2,n-2\alpha}$ if and only if 
$$
\|f\|_{{\Cal L}_{2,n-2\alpha}}=\sup_{I}\left(\big(\ell(I)\big)^{2\alpha-n}\int_{I}|f(x)-f_I|^2dx\right)^\frac12<\infty.
$$
\flushpar {\rm (ii)} $\dot{L}^2_\alpha=(-\Delta)^{-\frac\alpha{2}}L^2\hookrightarrow L^{\frac{2n}{n-2\alpha}}$ is sharp with 
$$
\sup_{\|f\|_{\dot{L}^2_\alpha}>0}
\frac{\|f\|_{L^{\frac{2n}{n-2\alpha}}}}{\|f\|_{\dot{L}^2_\alpha}}=\left(\frac{\Gamma\big(\frac{n-2\alpha}{2}\big)}{\Gamma\big(\frac{n+2\alpha}{2}\big)}\right)^\frac12
\left(\frac{\Gamma(n)}{\Gamma\big(\frac{n}{2}\big)}\right)^{\frac{\alpha}{n}}
\left(\int_{\rn}\frac{|e^{-2\pi y\cdot (1,0,...,0)}-1|^2}{|y|^{(n+2\alpha)}}dy\right)^{-\frac12}.
$$
\flushpar{\rm (iii)} $Q_{\alpha;\infty}^{-1}=\nabla\cdot(Q_\alpha)^n$, where a tempered distribution $f$ on $\rn$ belongs to $Q_{\alpha;\infty}^{-1}$ if and only if
$$
\|f\|_{Q^{-1}_{\alpha;\infty}}=\sup_{x\in\rn,r\in (0,\infty)}\left(r^{2\alpha-n}\int_0^{r^2}\int_{|y-x|<r}|e^{t\Delta}f(y)|^2t^{-\alpha}dydt\right)^\frac12<\infty.
$$
\endproclaim

Note that $\Cal L_{2,n-2\alpha}$ is the so-called Morrey space of square form (cf. [Ca] and [Pe] ) and $\Cal L_{2,n}=BMO$. So Theorem 1.1 (i) keeps true for $\alpha=0$ in the sense of $(-\Delta)^0 BMO=BMO$. Quite surprisingly, this part corresponds nicely to Strichartz's $(-\Delta)^{-\frac{\alpha}{2}}BMO$-equivalence [Str1, Theorem 3.3]: 
$$
f\in (-\Delta)^{-\frac{\alpha}{2}}BMO\Longleftrightarrow \sup_I \left(\big(\ell(I)\big)^{-n}\int_I\int_I\frac{|f(x)-f(y)|^2}{|x-y|^{n+2\alpha}}dxdy\right)^\frac12<\infty.
$$ 
The imbedding without best constant in Theorem 1.1 (ii) is well-known (see for example [MS, Theorem 1] and the related references therein) and very useful in the study of the semi-linear wave equations (cf. [LinS]). A close look at both (i) and (ii) reveals that $Q_\alpha$ behaves like an affine Sobolev space. In addition, Theorem 1.1 (iii) extends [KoTa, Theorem 1]: $BMO^{-1}=\nabla\cdot(BMO)^n$ that just says: $f\in BMO^{-1}$ if and only if there are $f_j\in BMO$ such that $f=\sum_{j=1}^n{\partial_j f_j}$. 

As with the latter, we recall that the Cauchy problem for the incompressible Navier-Stokes system on the half-space $\bn=(0,\infty)\times\rn$:
$$
\cases \partial_t u-\Delta u+(u\cdot\nabla)u-\nabla p=0,\quad\hbox{in}\quad \bn;\\
\nabla\cdot u=0,\quad \hbox{in}\quad \bn;\\
u|_{t=0}=a,\quad \hbox{in}\quad \rn
\endcases
\tag 1.4
$$
is to establish the existence of a solution (velocity) $u=u(t,x)=\big(u_1(t,x),...,u_n(t,x)\big)$ with a pressure $p=p(t,x)$ of the fluid at time $t\in [0,\infty)$ and position $x\in\rn$ that assumes the given data (initial velocity) $a=a(x)=(a_1(x),...,a_n(x))$. If the solution exists, is unique, and depends continuously on the initial data (with respect to a given topology), then we say that the Cauchy problem is well-posed in that topology.

Of particularly significant is the invariance of (1.4) under the scaling changes:
$$
\cases u(t,x)\mapsto u_\lambda(t,x)=\lambda u(\lambda^2t,\lambda x);\\
p(t,x)\mapsto p_\lambda(t,x)=\lambda^2p(\lambda^2t,\lambda x);\\
a(x)\mapsto a_\lambda(x)=\lambda a(\lambda x).
\endcases
$$
So if the trio $(u(t,x), p(t,x), a(x))$ satisfies (1.4) then the trio $(u_\lambda(t,x), p_\lambda(t,x), a_\lambda(x))$ is a solution of (1.4) for any $\lambda>0$. This leads to a consideration of the well-posedness for (1.4) with a Cauchy data being of the scaling invariance. Through the scale invariance 
$$
\|a_\lambda\|_{(L^n)^n}=\sum_{j=1}^n\|(a_\lambda)_j\|_{L^n}=\sum_{j=1}^n\|(a_j)_\lambda\|_{L^n}=\sum_{j=1}^n\|a_j\|_{L^n}=\|a\|_{(L^n)^n},
$$ 
Kato proved in [Ka] that (1.4) has mild solutions locally in time if $a\in (L^n)^n$ and globally if $\|a\|_{(L^n)^n}$ is small enough (for some generalizations of Kato's result, see e.g. [Ta] and [Y]). Furthermore, in [KoTa], Koch-Tataru found, among other results, that (1.4) still has mild solutions locally in time if $a\in (\overline{VMO^{-1}})^n$ and globally once
$$
\sum_{j=1}^n\|a_j\|_{BMO^{-1}}=\sum_{j=1}^n\sup_{x\in\rn,r\in (0,\infty)}\left(r^{-{n}}\int_0^{r^2}\int_{|y-x|<r}|e^{t\Delta}a_j(y)|^2dydt\right)^\frac12
$$  
is sufficiently small. Here and henceforward, by a mild solution $u(t,x)$ of (1.4) we mean that $u(t,x)$ solves the integral equation
$$
u(t,x)=e^{t\Delta}a(x)-\int_0^t e^{(t-s)\Delta}P\nabla\cdot(u\otimes u)ds,
$$ 
where $e^{t\Delta}a(x)=(e^{t\Delta}a_1(x),...,e^{t\Delta}a_n(x))$ and $P$ is the Helmboltz-Weyl projection: 
$$
P=\{P_{jk}\}_{j,k=1,...,n}=\{\delta_{jk}+R_jR_k\}_{j,k=1,...,n}
$$
with $\delta_{jk}$ being the Kronecker symbol and $R_j=\partial_j(-\Delta)^{-\frac12}$ being the Riesz transform. 

Observe that $\|\cdot\|_{BMO^{-1}}$ and $\|\cdot\|_{Q_{\alpha;\infty}^{-1}}$ are also invariant under the scale transform $a(x)\mapsto \lambda a(\lambda x)$. So it is a natural thing to extend the results of Kato and Koch-Tataru to the $Q_\alpha$-setting. To do this, we introduce the following concept. 

\proclaim{Definition 1.2} Let $\alpha\in (0,1)$ and $T\in (0,\infty]$. Then we say: 

\flushpar{\rm (i)} A tempered distribution $f$ on $\rn$ belongs to the space $Q_{\alpha;T}^{-1}$ provided
$$
\|f\|_{Q^{-1}_{\alpha;T}}=\sup_{x\in\rn,r\in (0,{T})}\left(r^{2\alpha-n}\int_0^{r^2}\int_{|y-x|<r}|e^{t\Delta}f(y)|^2t^{-\alpha}dydt\right)^\frac12<\infty;
$$

\flushpar{\rm (ii)} A tempered distribution $f$ on $\rn$ belongs to $\overline{VQ_\alpha^{-1}}$ provided $\lim_{T\to 0}\|f\|_{Q^{-1}_{\alpha;T}}=0$;

\flushpar{\rm (iii)} A function $g$ on $\bn$ belongs to the space $X_{\alpha; T}$ provided
$$
\align
&\|g\|_{X_{\alpha; T}}\\
&=\sup_{t\in (0,T)}\sqrt{t}\|g(t,\cdot)\|_{L^\infty}+\sup_{x\in\rn,r^2\in (0,{T})}
\left(r^{2\alpha-n}\int_0^{r^2}\int_{|y-x|<r}|g(t,y)|^2t^{-\alpha}dydt\right)^\frac12<\infty.
\endalign
$$
In particular, we write 
$$
Q^{-1}_{0;T}=BMO^{-1}_T,\quad\overline{VQ_{0}^{-1}}=\overline{VMO^{-1}}\quad\hbox{and}\quad X_{0;T}=X_T.
$$ 
\endproclaim

Two immediate comments are given below: If 
$$
f_\lambda(x)=\lambda f(\lambda x)\quad\hbox{and}\quad g_\lambda(t,x)=\lambda g(\lambda^2 t,\lambda x)\quad\hbox{for}\quad \lambda,t>0\quad\hbox{and}\quad x\in\rn,
$$
then 
$$
\|f_\lambda\|_{Q^{-1}_{\alpha; \infty}}=\|f\|_{Q^{-1}_{\alpha; \infty}}\ \ \hbox{and}\ \ \|g_\lambda\|_{X_{\alpha; \infty}}=\|g\|_{X_{\alpha; \infty}};
$$
that is, $\|\cdot\|_{Q^{-1}_{\alpha; \infty}}$ and $\|\cdot\|_{X_{\alpha; \infty}}$ are scaling invariant. Second, we have $$
L^n\subseteq Q^{-1}_{\alpha;1}\subseteq BMO^{-1}_1\ \ \hbox{and}\ \ X_{\alpha;1}\subseteq X_1.
$$ 
To see this, note that $\|f\|_{BMO^{-1}}=\|f\|_{Q^{-1}_{0;\infty}}\le\|f\|_{Q^{-1}_{\alpha;\infty}}.$ 
Additionally, recall that $f\in \dot{B}_{p,\infty}^{-1+\frac np},\ p>n$ if and only if
$\|e^{t\Delta}f\|_{L^p}\lesssim t^{\frac{n-p}{2p}}$ for $t\in (0,1].$
Using H\"older's inequality, we obtain that if $f\in \dot{B}_{p,\infty}^{-1+\frac np},\ p>n$ and $r\in (0,1)$, then 
$$
\int_0^{r^2}\int_{|y-x|<r}|e^{t\Delta}f(y)|^2t^{-\alpha}dydt\lesssim r^{n(p-2)}\int_0^{r^2}\|e^{t\Delta}f\|_{L^p}^2 t^{-\alpha}dt\lesssim r^{n-2\alpha}
$$
and hence $f\in Q^{-1}_{\alpha;1}$. This, together with the well-known inclusion (see e.g. [KoTa]): $L^n\subseteq \dot{B}_{p,\infty}^{-1+\frac np},\ p>n$, yields the desired inclusion.

Below is our result on the well-posedness for the incompressible Navier-Stokes system.

\proclaim{Theorem 1.2} Let $\alpha\in (0,1)$. Then
\smallskip
\flushpar {\rm (i)} The Navier-Stokes system (1.4) has a unique small global mild solution in $(X_{\alpha;\infty})^n$ for all initial data $a$ with $\nabla\cdot a=0$ and $\|a\|_{(Q^{-1}_{\alpha;\infty})^n}$ being small.
\smallskip
\flushpar{\rm (ii)} For any $T\in (0,\infty)$ there is an $\epsilon>0$ such that the Navier-Stokes system (1.4) has a unique small mild solution in $(X_{\alpha; T})^n$ on $(0,T)\times\rn$ when the initial data $a$ satisfies $\nabla\cdot a=0$ and $\|a\|_{(Q^{-1}_{\alpha;T})^n}\le\epsilon$. In particular for all $a\in \big(\overline{VQ_\alpha^{-1}}\big)^n$ with $\nabla\cdot a=0$ there exists a unique small local mild solution in $(X_{\alpha; T})^n$ on $(0,T)\times\rn$.
\endproclaim

In the case of $\alpha=0$, Theorem 1.2 goes back to Theorems 2-3 of Koch-Tataru in [KoTa]. However, it is perhaps worth pointing out that their Theorems 2-3 do not deduce our Theorem 1.2 even though $(Q^{-1}_{\alpha;T})^n$ and $\big(\overline{VQ_\alpha^{-1}}\big)^n$ are subspaces of $(BMO^{-1}_T)^n$ and $\big(\overline{VMO^{-1}}\big)^n$ respectively, since the $(X_{\alpha; T})^n$ is contained properly in $(X_{T})^n$ when $0<\alpha<1$.
\medskip

In the forthcoming two sections, we provide the proofs of the above-stated theorems. The argument of Theorem 1.1 (i) follows from a chain of integral estimates for singular integral operators (see e.g. [Ch] and [CoMS]) with being partially inspired by Wu-Xie's work [WuXi], while in the demonstration of Theorem 1.1 (ii) we formulate the integral involved in the Sobolev space as weighted integral of the Fourier transform of the given function and take Lieb's sharp estimate for the Hardy-Littlewood-Sobolev inequality into account. The justification of Theorem 1.1 (iii) is an extension of Koch-Tataru's argument for the $BMO$-setting in [KoTa]. In showing Theorem 1.2 (i)-(ii), we improve Lemari\'e-Rieusset's treatment (cf. [Le, Chapter 16]) of Koch-Tataru's proof of settling the case $\alpha=0$ (see again [KoTa, Theorems 2 and 3]) in order to handle any value $\alpha\in (0,1)$. More precisely, our proof rests on two technical lemmas of which Lemma 3.1 brings Schur's lemma into play and so makes a difference.

\heading
{2. Proof of Theorem 1.1}
\endheading

To verify Theorem 1.1 (i), we must understand the quadratic Morrey space in spirit of (1.2). 

\proclaim{Lemma 2.1} Given $\alpha\in (0,1)$. Let $f$ be a measurable complex-valued function on $\rn$. Then 
$$
f\in \Cal L_{2,n-2\alpha}\Longleftrightarrow\sup_{x\in\rn, r\in (0,\infty)}r^{2\alpha-n}\int_0^r\int_{|y-x|<r}{\Big|{\partial_t}e^{-t\sqrt{-\Delta}}f(y)\Big|^2}tdydt<\infty.\tag 2.1
$$ 
\endproclaim
\demo{Proof} Assume $f\in \Cal L_{2,n-2\alpha}$. Recall 
$$
\psi_0(x)=\frac{1+|x|^2-(n+1)\Gamma(\frac{n+1}{2})\pi^{-\frac{n+1}{2}}}{(1+|x|^2)^{\frac{n+3}{2}}}\quad\hbox{and}\quad (\psi_0)_t(x)=t\partial_t e^{-t\sqrt{-\Delta}}(x,0).
$$
For any ball $B=\{y\in\rn:\ |y-x|<r\}$ in $\rn$, let $2B$ be the double ball of $B$ and $f_{2B}=|2B|^{-1}\int_{2B}f$ be the mean value of $f$ on $2B$. Note that $|E|$ stands for the Lebesgue measure of a set $E$. Let also 
$$
f_1=(f-f_{2B})1_{2B},\quad f_2=(f-f_{2B})1_{\rn\setminus 2B}\quad\hbox{and}\quad f_3=f_{2B},
$$
where $1_E$ stands for the characteristic function of a set $E$. Since $\int_{\rn}\psi_0=0$, we conclude that
$$
t{\partial_t}e^{-t\sqrt{-\Delta}}f(y)=(\psi_0)_t\ast f(y)=(\psi_0)_t\ast f_1(y)+(\psi_0)_t\ast f_2(y).
$$
Concerning the first term $(\psi_0)_t\ast f_1(y)$, we have the following estimates
$$
\align
\int_0^r\int_B|(\psi_0)_t\ast f_1(y)|^2t^{-1}dydt
&\le\int_B\int_0^\infty|(\psi_0)_t\ast f_1(y)|^2t^{-1}dtdy\\
&\le\int_{\rn}\left(\int_0^\infty|(\psi_0)_t\ast f_1(y)|^2t^{-1}dt\right)dy\\
&\lesssim\int_{2B}|f(y)-f_{2B}|^2dy\\
&\lesssim r^{n-2\alpha}\|f\|_{{\Cal L}_{2,n-2\alpha}}^2,
\endalign
$$
where we have used the $L^2$-boundedness: $\|\Cal G(f_1)\|_{L^2}\lesssim\|f_1\|_{L^2}$ for the Littlewood-Paley $\Cal G$-function of $f_1$
$$
\Cal G(f_1)(y)=\left(\int_0^\infty\left|t{\partial_t}e^{-t\sqrt{-\Delta}}f_1(y)\right|^2t^{-1}dt\right)^\frac12.
$$

At the same time, if $(t,y)\in (0,r)\times B$ and $B_k$ is the ball with center $x$ and radius $2^kr$, then we take the Cauchy-Schwarz inequality into account, and obtain the following inequalities for $(\psi_0)_t\ast f_2(y)$:
$$
\align
|(\psi_0)_t\ast f_2(y)|
&\lesssim\int_{\rn\setminus 2B}\frac{t|f(z)-f_{2B}|}{(t+|x-z|)^{n+1}}dz\\
&\lesssim\int_{\rn\setminus 2B}\frac{t|f(z)-f_{2B}|}{(r+|x-z|)^{n+1}}dy\\
&\lesssim t\sum_{k=1}^\infty\int_{B_k}\frac{t|f(z)-f_{2B}|}{(r+|x-z|)^{n+1}}dz\\
&\lesssim t\sum_{k=1}^\infty(2^k r)^{-(n+1)}\int_{B_k}|f(z)-f_{2B}|dz\\
&\lesssim t r^{\alpha-n-1}\|f\|_{{\Cal L}_{2,n-2\alpha}}.
\endalign
$$
Consequently,
$$
\int_0^r\int_B|(\psi_0)_t\ast f_2(y)|^2t^{-1}dydt\lesssim r^{n-2\alpha}\|f\|_{{\Cal L}_{2,n-2\alpha}}^2.
$$
The above-established estimates yield that the supremum in (2.1) is finite.

To handle the converse part, denote by
$$
S(I)=\Big\{(t,x)\in\bn: \ t\in (0,\ell(I)]\quad\hbox{and}\quad x\in I\Big\}.
$$
the Carleson box based on a given cube $I\subset\rn$.

Suppose the supremum condition in (2.1) is satisfied. Then 
$$
\sup_{I}\big(\ell(I)\big)^{2\alpha-n}\int_{S(I)}{\big|{\partial_t}e^{-t\sqrt{-\Delta}}f(y)\big|^2}tdydt<\infty.\tag 2.2
$$ 
In order to verify $f\in \Cal L_{2,n-2\alpha}$, we consider the projection operator
$$
\Pi_{\psi_0} F(x)=\int_{\bn}F(t,y)(\psi_0)_t(x-y)t^{-1}dydt,
$$
and prove that if 
$$
\|F\|_{\Cal C_\alpha}=\sup_I \left(\big(\ell(I)\big)^{2\alpha-n}\int_{S(I)}|F(t,y)|^2t^{-1}dydt\right)^\frac12<\infty
$$
then for any cube $J\subset\rn$,
$$
\int_J|\Pi_{\psi_0} F(x)-\big(\Pi_{\psi_0} F\big)_J|^2dx\lesssim \big(\ell(J)\big)^{n-2\alpha}\|F\|_{\Cal C_\alpha}^2.\tag 2.3
$$
Given a cube $J\subset\rn$ and $\lambda>0$, define $\lambda J$ as the cube concentric with $J$ and with sidelength $\ell(\lambda J)=\lambda\ell(J)$. Let $F_1=F|_{S(2J)}$ and $F_2=F-F_1.$ Then by a result of Coifman-Mayer-Stein [CoMS, p. 328-329],
$$
\align
\int_J|\Pi_{\psi_0} F_1(x)|^2dx&\le\int_{\rn}|\Pi_{\psi_0} F_1(x)|^2dx\lesssim\int_{\bn}|F_1(t,y)|^2t^{-1}dydt\\
&\lesssim\int_{S(2J)}|F(t,y)|^2t^{-1}dydt\lesssim\|F\|_{\Cal C_\alpha}^2\big(\ell(J)\big)^{n-2\alpha}.
\endalign
$$
On the other hand, from the definition of $\Pi_{\psi_0}$, the boundedness of $\psi_0$ and the Cauchy-Schwarz inequality it follows that
$$
\align
&\int_J|\Pi_{\psi_0} F_2(x)|^2dx\\
&=\int_J\left|\int_{\bn}(\psi_0)_t(x-y)F_2(t,y)t^{-1}dydt\right|^2dx\\
&\le\int_J\left(\int_{\bn\setminus S(2J)}|(\psi_0)_t(x-y)||F_2(t,y)|t^{-1}dydt\right)^2dx\\
&\le\int_J\left(\sum_{k=1}^\infty\int_{S(2^{k+1}J)\setminus S(2^kJ)}|(\psi_0)_t(x-y)||F_2(t,y)|t^{-1}dydt\right)^2dx\\
&\lesssim\int_J\left(\sum_{k=1}^\infty\big(2^{k}\ell(J)\big)^{-1}\int_{S(2^{k+1}J)\setminus S(2^kJ)}|(\psi_0)_t(x-y)||F_2(t,y)|dydt\right)^2dx\\
&\lesssim\int_J\left(\sum_{k=1}^\infty\big(2^{k}\ell(J)\big)^{-(n+1)}\int_{S(2^{k+1}J)\setminus S(2^kJ)}|F_2(t,y)|dydt\right)^2dx\\
&\lesssim\int_J\left(\sum_{k=1}^\infty\big(2^{k}\ell(J)\big)^{-\frac{n}{2}}\left(\int_{S(2^{k+1}J)\setminus S(2^kJ)}|F_2(t,y)|^2t^{-1}dydt\right)^\frac12\right)^2dx\\
&\lesssim\|F\|_{\Cal C_\alpha}^2\big(\ell(J)\big)^{n-2\alpha}.
\endalign
$$
Thus
$$
\align
\int_J\big|\Pi_{\psi_0}F(x)-\big(\Pi_{\psi_0}F\big)_J\big|^2dx
&\lesssim\int_J\big|\Pi_{\psi_0}F(x)\big|^2dx\\
&\lesssim\int_J\big|\Pi_{\psi_0}F_1(x)\big|^2dx+\int_J\big|\Pi_{\psi_0}F_2(x)\big|^2dx\\
&\lesssim\|F\|_{\Cal C_\alpha}^2\big(\ell(J)\big)^{n-2\alpha},
\endalign
$$
namely, (2.3) holds. Applying (2.3) to $\Pi_{\psi_0}\big((\psi_0)_t\ast f\big)$ which equals $f$ under (2.2), we achieve $f\in \Cal L_{2,n-2\alpha}$.
\enddemo

\demo{\rm\bf Proof of Theorem 1.1 (i)} Since the imbedding with that constant can be readily derived from a routine computation, it suffices to show $Q_\alpha=(-\Delta)^{-\frac\alpha{2}}\Cal L_{2,n-2\alpha}$. 

For $f\in \Cal L_{2,n-2\alpha}$, let $F(t,y)=t^{1+\alpha}{\partial_t}e^{-t\sqrt{-\Delta}}f(y).$ 
Then by Lemma 2.1 we get that
$$
\sup_{x\in\rn, r\in (0,\infty)}r^{2\alpha-n}\int_0^r\int_{|y-x|<r}|F(t,y)|^2t^{-1-2\alpha}dydt\lesssim\|f\|_{\Cal L_{2,n-2\alpha}}^2,
$$
and so that $\Pi_{\psi_0}F\in Q_\alpha$ thanks to [DX1, Theorem 7.0 (i)]. Note however that
$$
\widehat{F(t,\cdot)}(x)=-t^{\alpha+2}|x|\hat{f}(x)\exp(-t|x|).
$$
So a calculation infers
$$
\widehat{\Pi_{\psi_0}F}(x)=2^{-(2+\alpha)}\Gamma(2+\alpha)|x|^{-\alpha}\hat{f}(x)={2^{-2}\pi^\alpha\Gamma(2+\alpha)}\widehat{(-\Delta)^{-\frac{\alpha}{2}}f}(x).
$$
Therefore, $(-\Delta)^{-\frac{\alpha}{2}}f$ belongs to $Q_\alpha$.

Conversely, suppose $g\in Q_\alpha$. Setting $G(t,y)=t^{1-\alpha}{\partial_t}e^{-t\sqrt{-\Delta}}g(y)$, we deduce $\|G\|_{\Cal C_\alpha}\lesssim\|g\|_{Q_\alpha}$ by using (1.2). Thus (2.3) is valid for this $G(\cdot,\cdot)$. From the argument of Lemma 2.1 it is easily derived that $\Pi_{\psi_0}G\in \Cal L_{2,n-2\alpha}$. Since
$$
\widehat{\Pi_{\psi_0}G}(x)={2^{-2}\pi^{-\alpha}}{\Gamma(2-\alpha)}\widehat{(-\Delta)^{\frac{\alpha}{2}}g}(x),
$$
we conclude that $f=(-\Delta)^{\frac{\alpha}{2}}g\in \Cal L_{2,n-2\alpha}$ and $g=(-\Delta)^{-\frac{\alpha}{2}}f.$
This completes the proof.
\enddemo

\demo{\rm\bf Proof of Theorem 1.1 (ii)} According to [Str2, p. 175], we use Fubini's theorem, Plancherel's formula, the change of variables $y=|x|^{-1}z$ and an orthonormal transform to obtain
$$
\align
\|f\|_{\dot{L}^2_\alpha}^2&=\int_{\rn}\left(\int_{\rn}|f(x+y)-f(x)|^2dx\right)|y|^{-(n+2\alpha)}dy\\
&=\int_{\rn}\left(\int_{\rn}|e^{-2\pi y\cdot x}-1|^2|y|^{-(n+2\alpha)}dy\right)|\hat{f}(x)|^2dx\\
&=\left(\int_{\rn}|e^{-2\pi z\cdot (1,0,...,0)}-1|^2|z|^{-(n+2\alpha)}dz\right)\int_{\rn}|\hat{f}(x)|^2|x|^{2\alpha}dx.
\endalign
$$
Accordingly, $\dot{L}^2_\alpha=(-\Delta)^{-\frac\alpha{2}}L^2$. Note that $\hat{f}(x)=\int_{\rn}f(y)\exp(-2\pi i x\cdot y)dy$ and 
$$
e^{-t\sqrt{-\Delta}}f(x)=\int_{\rn}\hat{f}(y)\exp\big(-2\pi(iy\cdot x+|y|t)\big)dy.
$$
So, by differentiation and integration (cf. [Ste, p. 83]), we find
$$
\|\nabla e^{-t\sqrt{-\Delta}}f\|_{L^2}^2=8\pi^2\int_{\rn}|x|^2|\hat{f}(x)|^2\exp(-4\pi|x|t)dx.
$$
Consequently,
$$
\int_0^\infty\|\nabla e^{-t\sqrt{-\Delta}}f\|^2_{L^2}t^{1-2\alpha}dt=\frac{8\pi^2\Gamma\big(2(1-\alpha)\big)}{(4\pi)^{2(1-\alpha)}}\int_{\rn}|x|^{2\alpha}|\hat{f}(x)|^2dx.\tag 2.4
$$
Moreover, the preceding consideration actually tells us that proving the desired sharp imbedding amounts to verifying the following best inequality 
$$
\|f\|_{L^{\frac{2n}{n-2\alpha}}}^2\le \tau_{n,\alpha}\int_0^\infty\|\nabla e^{-t\sqrt{-\Delta}}f\|^2_{L^2}t^{1-2\alpha}dt,\tag 2.5
$$
where
$$
\tau_{n,\alpha}=\frac{2^{1-4\alpha}\Gamma\Big(\frac{n-2\alpha}{2}\Big)}{\pi^\alpha\Gamma\big(2(1-\alpha)\big)\Gamma\Big(\frac{n+2\alpha}{2}\Big)}\left(\frac{\Gamma(n)}{\Gamma\big(\frac{n}2\big)}\right)^\frac{2\alpha}{n}.
$$
To this end, we use $\langle f,g\rangle$ as $\int_{\rn}f(x)\overline{g(x)}dx$, and then get 
$$
|\langle f,g\rangle|^2=|\langle \hat{f},\hat{g}\rangle|^2\le\left(\int_{\rn}|x|^{2\alpha}|\hat{f}(x)|^2dx\right) \left(\int_{\rn}|x|^{-2\alpha}|\hat{g}(x)|^2dx\right).
$$
Because (cf. [LieL, Corollary 5.10])
$$
\int_{\rn}|x|^{-2\alpha}|\hat{g}(x)|^2dx=\frac{\pi^{2\alpha-\frac{n}{2}}\Gamma\big(\frac{n-2\alpha}{2}\big)}{\Gamma(\alpha)}\int_{\rn}\int_{\rn}\frac{g(x)\overline{g(y)}}{|x-y|^{n-2\alpha}}dxdy,
$$
it follows from Lieb's sharp version [Lie] of the Hardy-Littlewood-Sobolev inequality that
$$
|\langle f,g\rangle|^2\le\frac{\pi^{{\alpha}}{\Gamma\big(\frac{n-2\alpha}{2}\big)}}{\Gamma\big(\frac{n+2\alpha}{2}\big)}\left(\frac{\Gamma(n)}{\Gamma\big(\frac{n}{2}\big)}\right)^\frac{2\alpha}n \|g\|^2_{L^{\frac{2n}{n+2\alpha}}}\int_{\rn}|x|^{2\alpha}|\hat{f}(x)|^2dx.
$$
In the last inequality we take $g=f|f|^{\frac{4\alpha}{n-2\alpha}}$ for $f\in \dot{L}^2_\alpha$, and use (2.4) to achieve (2.5) whose equality can be checked for $f(x)=(1+|x|^2)^{(2\alpha-n)/2}$ through a direct calculation.
\enddemo

The proof of Theorem 1.1 (iii) depends on the following lemma.

\proclaim{Lemma 2.2} Given $\alpha\in (0,1)$. For $j,k=1,...,n$ let $f_{j,k}={\partial_j\partial_k} (-\Delta)^{-1}f$. If $f\in Q_{\alpha;\infty}^{-1}$ then $f_{j,k}\in Q^{-1}_{\alpha;\infty}$.
\endproclaim
\demo{Proof} Assume that $\phi\in\Cal D(\rn)$ is fixed with $\hbox{supp}\phi\subset B(0,1)=\{x\in\rn: |x|<1\}$ and $\int_{\rn}\phi=1$. Recall $\phi_r(x)=r^{-n}\phi(x/r)$, and write $g_r(t,x)=\phi_r\ast\partial_j\partial_k (-\Delta)^{-1}e^{t\Delta}f(x)$. Then
$$
e^{t\Delta}f_{j,k}(x)=\partial_j\partial_k (-\Delta)^{-1}e^{t\Delta}f(x)=f_r(t,x)+g_r(t,x).
$$
If $\dot{B}^{1,1}_1$ stands for the predual of the homogeneous Besov space $\dot{B}^{-1,\infty}_\infty$, then $f\in Q^{-1}_{\alpha;\infty}$ yields $f\in BMO^{-1}\subseteq\dot{B}^{-1,\infty}_\infty$ (see also [Le, p. 160, Lemma 16.1]) and
$$
\|g_r(t,\cdot)\|_{L^\infty}\le\|\phi_r\|_{\dot{B}^{1,1}_1}\Big\|\partial_j\partial_k (-\Delta)^{-1}e^{t\Delta}f\Big\|_{\dot{B}^{-1,\infty}_\infty}\lesssim r^{-1}\|f\|_{\dot{B}^{-1,\infty}_\infty}.
$$
Consequently, 
$$
\int_0^{r^2}\int_{|y-x|<r}|g_r(t,y)|^2t^{-\alpha}dydt\lesssim r^{n-2\alpha}\|f\|^2_{\dot{B}^{-1,\infty}_\infty}\lesssim r^{n-2\alpha}\|f\|^2_{Q^{-1}_{\alpha;\infty}}.\tag 2.6
$$
Next, we estimate $f_r$. Take $\psi\in\Cal D(\rn)$ so that $\psi=1$ on $B(0,10)=\{x\in\rn: |x|<10\}$. Define $\psi_{r,x}=\psi\big(\frac{y-x}{r}\big)$ and write $f_r=F_{r,x}+G_{r,x}$ where
$$
G_{r,x}={\partial_j\partial_k}(-\Delta)^{-1}\psi_{r,x}e^{t\Delta}f-\phi_r\ast{\partial_j\partial_k}(-\Delta)^{-1}\psi_{r,x}e^{t\Delta}f.
$$
Thus, we employ the Plancherel formula for the space variable to get
$$
\align
&\int_0^{r^2}\big\|{\partial_j\partial_k}(-\Delta)^{-1}\psi_{r,x}e^{t\Delta}f\big\|_{L^2}^2\frac{dt}{t^\alpha}\\
&\lesssim\int_0^{r^2}\big\|\widehat{(\partial_j\partial_k(-\Delta)^{-1}\psi_{r,x}e^{t\Delta}f\big)}\big\|_{L^2}^2\frac{dt}{t^\alpha}\\
&\lesssim\int_0^{r^2}\int_{\rn}\big|y_jy_k|y|^{-2}\widehat{\big(\psi_{r,x}e^{t\Delta}f\big)}(y)\big|^2t^{-\alpha}dydt\\
&\lesssim\int_0^{r^2}\big\|\widehat{\psi_{r,x}e^{t\Delta}f}\big\|_{L^2}^2\frac{dt}{t^\alpha}\\
&\lesssim\int_0^{r^2}\big\|{\psi_{r,x}e^{t\Delta}f}\big\|_{L^2}^2\frac{dt}{t^\alpha}.
\endalign
$$
And, by Minkowski's inequality (for $\phi_r$) as well as the Plancherel formula again, we have
$$
\int_0^{r^2}\big\|\phi_r\ast{\partial_j\partial_k} (-\Delta)^{-1}\psi_{r,x}e^{t\Delta}f\big\|_{L^2}^2\frac{dt}{t^\alpha}\lesssim\int_0^{r^2}\big\|{\psi_{r,x}e^{t\Delta}f}\|_{L^2}^2\frac{dt}{t^\alpha}.
$$
The last two estimates imply
$$
\int_0^{r^2}\big\|G_{r,\cdot}(t,\cdot)\big\|_{L^2}^2\frac{dt}{t^\alpha}\lesssim\int_0^{r^2}\big\|{\psi_{r,x}e^{t\Delta}f}\big\|_{L^2}^2\frac{dt}{t^\alpha}.
$$
To control $F_{r,x}$, we bring the following estimate proved in [Le, p. 161] 
$$
\int_{|y-x|<r}|F_{r,x}(t,y)|^2dy\lesssim r^{n+1}\int_{|w-x|\ge 10r}{|e^{t\Delta}f(w)|^2}{|x-w|^{-(n+1)}}dw
$$
into play, and get
$$
\align
&\int_0^{r^2}\int_{|y-x|<r}|F_{r,x}(t,y)|^2t^{-\alpha}dydt\\
&\lesssim r^{n+1}\int_{|w-x|\ge 10r}{|x-w|^{-(n+1)}}\int_0^{r^2}|e^{t\Delta}f(w)|^2t^{-\alpha}dtdw\\
&\lesssim  r^{n+1}\sum_{l=1}^\infty\int_{10^lr\le|w-x|\le 10^{l+1}r}{|x-w|^{-(n+1)}}\int_0^{r^2}|e^{t\Delta}f(w)|^2t^{-\alpha}dtdw\\
&\lesssim \sum_{l=1}^\infty{10^{-l(n+1)}}\int_{|w-x|\le 10^{l+1}r}\int_0^{r^2}|e^{t\Delta}f(w)|^2t^{-\alpha}dtdw\\
&\lesssim r^{n-2\alpha}\|f\|^2_{Q^{-1}_{\alpha;\infty}}.
\endalign
$$
The integral estimates on $F_{r,x}$ and $G_{r,x}$ give
$$
\int_0^{r^2}\int_{|y-x|<r}|f_{r}(t,y)|^2t^{-\alpha}dydt\lesssim r^{n-2\alpha}\|f\|^2_{Q^{-1}_{\alpha;\infty}}.\tag 2.7
$$
Combining (2.6) and (2.7) gives $f_{j,k}\in Q_{\alpha;\infty}^{-1}$, as required.
\enddemo

\demo{\rm\bf Proof of Theorem 1.1 (iii)} If $f\in\nabla\cdot(Q_{\alpha})^n$, then there are $f_1,...,f_n\in Q_\alpha$ such that
$f=\sum_{j=1}^n{\partial_j}f_j$. So, the Minkowski inequality derives
$$
\|f\|_{Q_{\alpha;\infty}^{-1}}\le\sum_{j=1}^n\big\|{\partial_j}f_j\big\|_{Q_{\alpha;\infty}^{-1}}\lesssim\sum_{j=1}^n\|f_j\|_{Q_{\alpha}}.
$$
This means $f\in Q_{\alpha;\infty}^{-1}$. 

Conversely, let $f\in Q_{\alpha;\infty}^{-1}$. If $f_{j,k}={\partial_j\partial_k}(-\Delta)^{-1}f,\quad j,k=1,...,n,$ then $f_{j,k}\in Q^{-1}_{\alpha;\infty}$ by Lemma 2.2, and hence $f_k=-{\partial_k}(-\Delta)^{-1}f\in Q_\alpha$. This leads to 
$$
\widehat{\sum_{k=1}^n{\partial_k}f_k}=-\sum_{k=1}^n\widehat{f_{k,k}}=\hat{f},
$$
completing the proof.
\enddemo

\heading
{3. Proof of Theorem 1.2}
\endheading

To prove Theorem 1.2 we need two lemmas.

\proclaim{Lemma 3.1} Let $\alpha\in (0,1)$. Given a number $T\in (0,\infty]$ and a function $f(\cdot,\cdot)$ on $\bn$. Let $Af(t,x)=\int_0^t e^{(t-s)\Delta}\Delta f(s,x)ds$. Then 
$$
\int_0^T\big\|Af(t,\cdot)\big\|_{L^2}^2\frac{dt}{t^\alpha}\lesssim\int_0^T\big\|f(t,\cdot)\big\|_{L^2}^2\frac{dt}{t^\alpha}.\tag 3.1
$$
\endproclaim
\demo{Proof} It suffices to justify (3.1) for $T=\infty$. This is because: If $T<\infty$, then one may extend $f$ by putting $f=0$ on $(T,\infty)$, since $Af$ counts only on the values of $f$ on $(0,t)\times\rn$. Moreover, we may define $f=0=Af$ for $t\in (-\infty,0)$. 

Recall $e^{t\Delta}(x,0)=(4\pi t)^{-\frac n2}\exp(-\frac{|x|^2}{4t})$. Define
$$
\Omega(t,x)=\cases \Delta e^{t\Delta}(x,0),\quad & t>0\\
0,\quad & t\le 0.
\endcases
$$
Then we read that
$$
Af(t,x)=\int_{\Bbb R}\int_{\rn}\Omega(t-s,x-y)f(s,y)dyds,
$$
and hence $A$ becomes a convolution operator over $\Bbb R^{1+n}=\Bbb R\times\rn$. Since
$$
\widehat{\Omega(t,\cdot)}(\zeta)=\int_{\rn}\Omega(t,x)\exp(-2\pi i x\cdot \zeta)dx=-(2\pi)^2|\zeta|^2\exp\big(-(2\pi)^2t|\zeta|^2\big),\quad t>0,
$$
we conclude that
$$
\align
\widehat{Af(t,\cdot)}(\zeta)&=\int_{\Bbb R^{1+n}}\left(\int_{\rn}\Omega(t-s,x-y)f(s,y)\exp(-2\pi i x\cdot\zeta)dx\right)dyds\\
&=\int_{\Bbb R^{1+n}}f(s,y)\left(\int_{\rn}\Omega(t-s,u)\exp(-2\pi i (u+y)\cdot\zeta)du\right)dyds\\
&=\int_0^t\int_{\rn}f(s,y)\exp(-2\pi i y\cdot\zeta)\Big(-(2\pi)^2|\zeta|^2\exp\big(-(2\pi)^2(t-s)|\zeta|^2\big)\Big)dyds\\
&=-(2\pi)^2\int_0^t|\zeta|^2\exp\big(-(2\pi)^2(t-s)|\zeta|^2\big)\widehat{f(s,\cdot)}(\zeta)ds.
\endalign
$$
This formula, together with the Fubini theorem and the Plancherel formula, implies 
$$
\align
&\int_0^\infty\big\|Af(t,\cdot)\big\|^2_{L^2}\frac{dt}{t^\alpha}\\
&=\int_0^\infty\left(\int_{\rn}|\widehat{Af(t,\cdot)}(x)|^2dx\right)t^{-\alpha}dt\\
&\le\int_0^\infty\left(\int_{\rn}\Big(\int_0^t\frac{|\zeta|^2}{\exp\big((2\pi)^2(t-s)|\zeta|^2\big)}|\widehat{f(s,\cdot)}(\zeta)|ds\Big)^2d\zeta\right)t^{-\alpha}dt\\
&=(2\pi)^2\int_{\rn}\left(\int_0^\infty\Big(\int_0^t\frac{|\zeta|^2}{\exp\big((2\pi)^2(t-s)|\zeta|^2\big)}|\widehat{f(s,\cdot)}(\zeta)|ds\Big)^2t^{-\alpha}dt\right)d\zeta\\
&=(2\pi)^2\int_{\rn}\left(\int_0^\infty\Big(\int_0^\infty \big(1_{\{0\le s\le t\}}\big)\frac{|\zeta|^2}{\exp\big((2\pi)^2(t-s)|\zeta|^2\big)}|\widehat{f(s,\cdot)}(\zeta)|ds\Big)^2t^{-\alpha}dt\right)d\zeta.
\endalign
$$
This tells us that if one can prove
$$
\int_0^\infty\left(\int_0^\infty \Big(1_{\{0\le s\le t\}}\Big)\frac{|\zeta|^2|\widehat{f(s,\cdot)}(\zeta)|}{\exp\big((t-s)|\zeta|^2\big)}ds\right)^2t^{-\alpha}dt\lesssim\int_0^\infty |\widehat{f(t,\cdot)}(\zeta)|^2t^{-\alpha}dt,\tag 3.2
$$
then the Plancherel formula can be used again to yield
$$
\int_0^\infty\big\|Af(t,\cdot)\big\|_{L^2}^2\frac{dt}{t^\alpha}\lesssim\int_0^\infty\big\|f(t,\cdot)\|_{L^2}^2\frac{dt}{t^\alpha},
$$
as desired. 

To verify (3.2), we rewrite its left side as
$$
\int_0^\infty\left(\int_0^\infty K(s,t)F(s,\zeta)ds\right)^2dt,
$$
where
$$
F(s,\zeta)=s^{-\frac{\alpha}{2}}|\widehat{f(s,\cdot)}(\zeta)|\quad\hbox{and}\quad K(s,t)=\big(1_{\{0\le s\le t\}}\big)\Big(\frac{s}{t}\Big)^{\frac{\alpha}{2}}\frac{|\zeta|^2}{\exp((t-s)|\zeta|^2)}.
$$
Clearly,
$$
\int_0^\infty K(s,t)ds=|\zeta|^2\int_0^t\Big(\frac{s}{t}\Big)^{\frac{\alpha}{2}}\exp(-(t-s)|\zeta|^2)ds\le 1-\exp(-t|\zeta|^2)\le 1
$$
and
$$
\int_0^\infty K(s,t)dt=|\zeta|^2\int_s^\infty\Big(\frac{s}{t}\Big)^{\frac{\alpha}{2}}\exp(-(t-s)|\zeta|^2)dt\le 1.
$$
So, by Schur's lemma we get 
$$
\int_0^\infty\left(\int_0^\infty K(s,t)F(s,\zeta)ds\right)^2dt\lesssim\int_0^\infty \big(F(t,\zeta)\big)^2dt,
$$
reaching (3.2).
\enddemo

\proclaim{Lemma 3.2} Given $\alpha\in (0,1)$. For a function $f$ on $(0,1)\times\rn$ let
$$
C(f;\alpha)=\sup_{x\in\rn, r\in (0,1)}r^{2\alpha-n}\int_0^{r^2}\int_{|y-x|<r}|f(t,y)|t^{-\alpha}dtdy.
$$
Then
$$
\int_0^1\left\|\sqrt{-\Delta}e^{t\Delta}\int_0^t f(s,\cdot)ds\right\|_{L^2}^2\frac{dt}{t^\alpha}\lesssim C(f;\alpha)\int_0^1\big\|f(t,\cdot)\big\|_{L^2}^2\frac{dt}{t^\alpha}.\tag 3.3
$$
\endproclaim
\demo{Proof} In the sequel, $\langle\cdot,\cdot\rangle$ stands for the inner product in $L^2$ with respect to the space variable $x\in\rn$. So
$$
\align
\|\cdots\|_{L^2}^2&=\int_{\rn}\left|\sqrt{-\Delta}e^{t\Delta}\int_0^t f(s,y)ds\right|^2 dy\\
&=\left\langle \sqrt{-\Delta}e^{t\Delta}\int_0^t f(s,y)ds, \sqrt{-\Delta}e^{t\Delta}\int_0^t f(s,y)ds\right\rangle\\
&=\int_0^t\int_0^t\left\langle \sqrt{-\Delta}e^{t\Delta}f(s,\cdot), \sqrt{-\Delta}e^{t\Delta}f(h,\cdot)\right\rangle dsdh.
\endalign
$$
This gives
$$
\align
\int_0^1\|\cdots\|_{L^2}^2\frac{dt}{t^\alpha}&=2\Re\left(\iint_{0<h<s<1}\left\langle f(s,\cdot), \int_s^1(-\Delta)e^{2t\Delta}f(h,\cdot)t^{-\alpha}dt\right\rangle dsdh\right)\\
&\lesssim\iint_{0<h<s<1}\left\langle |f(s,\cdot)|, (e^{2\Delta}-e^{2s\Delta})|f(h,\cdot)|\right\rangle dsdh\\
&\lesssim\int_0^1\left\langle |f(s,\cdot)|, \int_0^s(e^{2\Delta}-e^{2s\Delta})|f(h,\cdot)|dh\right\rangle s^{-\alpha}ds\\
&\lesssim\left(\int_0^1\|f(s,\cdot)\|_{L^1}s^{-\alpha}ds\right)\sup_{s\in (0,1]}\left\|\int_0^s e^{2s\Delta}|f(h,\cdot)|dh\right\|_{L^\infty}.
\endalign
$$
As estimated in [Le, p. 163], it follows that 
$$
\sup_{z\in\rn, s\in (0,1]}\int_0^s e^{2s\Delta}|f(h,z)|dh
\lesssim\sup_{x\in\rn, r\in (0,1)}r^{-n}\int_0^{r^2}\int_{|y-x|<r}|f(s,y)|dyds.
$$
This estimate in turn implies 
$$
\sup_{z\in\rn, s\in (0,1]}\int_0^s e^{2s\Delta}|f(h,z)|dh\lesssim\sup_{x\in\rn, r\in (0,1)}r^{2\alpha-n}\int_0^{r^2}\int_{|y-x|<r}|f(s,y)|s^{-\alpha}dsdy
$$
and 
$$
\int_0^1\|\cdots\|_{L^2}^2\frac{dt}{t^\alpha}\lesssim 
C(f;\alpha)\left(\int_0^1\|f(s,\cdot)\|_{L^1}s^{-\alpha}ds\right),
$$
giving (3.3).
\enddemo

\demo{\rm\bf Proof of Theorem 1.2 (i)-(ii)} In accordance with the Picard contraction principle (cf. [Le, 145, Theorem 15.1]), we find that proving Theorem 1.2 amounts to demonstrating that the bilinear operator
$$
B(u,v)=\int_0^t e^{(t-s)\Delta} P\nabla\cdot(u\otimes v)ds
$$
is bounded from $(X_{\alpha;T})^n\times(X_{\alpha;T})^n$ to $(X_{\alpha;T})^n$. Naturally, $u\in (X_{\alpha;T})^n$ and $a\in (Q^{-1}_{\alpha;T})^n$ are respectively equipped with the norms
$$
\|u\|_{(X_{\alpha;T})^n}=\sum_{j=1}^n\|u_j\|_{X_{\alpha;T}}\quad\hbox{and}\quad\|a\|_{(Q^{-1}_{\alpha;T})^n}=\sum_{j=1}^n\|a_j\|_{Q^{-1}_{\alpha;T}}.
$$ 

{\it Step 1. $L^\infty$-bound}. We are about to prove that if $t\in (0,T)$ then
$$
|B(u,v)|\lesssim t^{-\frac12}\|u\|_{(X_{\alpha;T})^n}\|v\|_{(X_{\alpha;T})^n}.\tag 3.4
$$
If $\frac{t}{2}\le s<t$ then
$$
\|e^{(t-s)\Delta} P\nabla\cdot(u\otimes v)\|_{L^\infty}\lesssim (t-s)^{-\frac12}\|u\|_{L^\infty}\|v\|_{L^\infty}\lesssim (t-s)^{-\frac12}s^{-1}\|u\|_{(X_{\alpha;T})^n}\|v\|_{(X_{\alpha;T})^n}.
$$
If $0<s<\frac{t}{2}$ then 
$$
\align
&|e^{(t-s)\Delta}P\nabla\cdot(u\otimes v)|\\
&\lesssim\int_{\rn}\frac{|u(s,y)||v(s,y)|}{(\sqrt{t}+|x-y|)^{n+1}}dy\\
&\lesssim\sum_{{k}\in\Bbb Z^n}(\sqrt{t}(1+|k|)^{-(n+1)}\int_{x-y\in\sqrt{t}(k+[0,1]^n)}|u(s,y)||v(s,y)|dy.
\endalign
$$
An application of the Cauchy-Schwarz inequality yields
$$
\align
&\int_0^t\int_{x-y\in\sqrt{t}(k+[0,1]^n)}|u(s,y)||v(s,y)|dyds\\
&\lesssim t^\alpha\left(\int_0^t\int_{x-y\in\sqrt{t}(k+[0,1]^n)}\frac{|u(s,y)|^2}{s^{\alpha}}{dyds}\right)^\frac12\left(\int_0^t\int_{x-y\in\sqrt{t}(k+[0,1]^n)}\frac{|v(s,y)|^2}{s^{\alpha}}dyds\right)^\frac12\\
&\lesssim t^\frac{n}{2}\|u\|_{(X_{\alpha;T})^n}\|v\|_{(X_{\alpha;T})^n}.
\endalign
$$
From the foregoing inequalities it follows that
$$
\align
|B(u,v)|&\lesssim\int_0^{\frac t2}|e^{(t-s)\Delta}P\nabla\cdot(u\otimes v)|ds+\int_{\frac t2}^t|e^{(t-s)\Delta}P\nabla\cdot(u\otimes v)|ds\\
&\lesssim t^{-\frac12}\|u\|_{(X_{\alpha;T})^n}\|v\|_{(X_{\alpha;T})^n}+\left(\int_{\frac{t}{2}}^t s^{-1}(t-s)^{-\frac12}ds\right)\|u\|_{(X_{\alpha;T})^n}\|v\|_{(X_{\alpha;T})^n}\\
&\lesssim t^{-\frac12}\|u\|_{(X_{\alpha;T})^n}\|v\|_{(X_{\alpha;T})^n},
\endalign
$$
establishing (3.4).

{\it Step 2. $L^2$-bound}. We are about to show that if $x\in\rn$ and $r^2\in (0,T)$ then
$$
r^{2\alpha-n}\int_0^{r^2}\int_{|y-x|<r}|B(u,v)|^2s^{-\alpha}dyds\lesssim\|u\|^2_{(X_{\alpha;T})^n}\|v\|^2_{(X_{\alpha;T})^n}.\tag 3.5
$$
To do so, let $1_{r,x}(y)=1_{\{|y-x|<10 r\}}(y)$, i.e., the characteristic function on the ball $\{y\in\rn: |y-x|<10 r\}$, and set $B(u,v)=B_1(u,v)-B_2(u,v)-B_3(u,v)$, where
$$
B_1(u,v)=\int_0^s e^{(s-h)\Delta}P\nabla\cdot\big((1-1_{r,x})u\otimes v\big)dh,
$$
$$
B_2(u,v)=(-\Delta)^{-\frac 12}P\nabla\cdot\int_0^s e^{(s-h)\Delta}\Delta \Big((-\Delta)^{-\frac12}(I-e^{h\Delta})(1_{r,x})u\otimes v\Big)dh,
$$
and
$$
B_3(u,v)=(-\Delta)^{-\frac 12}P\nabla\cdot(-\Delta)^\frac12 e^{s\Delta}\Big(\int_0^s\big(1_{r,x})u\otimes v\big)dh\Big).
$$
Here and henceafter, $I$ stands for the identity operator.

When $0<s<r^2$ and $|y-x|<r$, the Cauchy-Schwarz inequality produces
$$
\align
&|B_1(u,v)|\\
&\lesssim\int_0^s\int_{|z-x|\ge 10 r}\frac{|u(h,z)||v(h,z)|}{(\sqrt{s-h}+|y-z|)^{n+1}}dzdh\\
&\lesssim\int_0^{r^2}\int_{|z-x|\ge 10 r}\frac{|u(h,z)||v(h,z)|}{|x-z|^{n+1}}dzdh\\
&\lesssim\left(\int_0^{r^2}\int_{|z-x|\ge 10 r}\frac{|u(h,z)|^2}{|x-z|^{n+1}}dzdh\right)^\frac12\left(\int_0^{r^2}\int_{|y-x|\ge 10 r}\frac{|v(h,z)|^2}{|x-z|^{n+1}}dzdh\right)^\frac12\\
&\lesssim r^{-1}\|u\|_{(X_{\alpha;T})^n}\|v\|_{(X_{\alpha;T})^n}.
\endalign
$$
Therefore,
$$
\int_0^{r^2}\int_{|y-x|<r}|B_1(u,v)|^2t^{-\alpha}dydt\lesssim r^{n-2\alpha}\|u\|_{(X_{\alpha;T})^n}^2\|v\|_{(X_{\alpha;T})^n}^2.
$$

For $B_2(u,v)$, put
$$
M(h,y)=1_{r,x}(u\otimes v)=1_{r,x}(y)\big(u(h,y)\otimes v(h,y)\big).
$$
By the $L^2$-boundedness of the Riesz transform and Lemma 3.1 we achieve 
$$
\align
\int_0^{r^2}\big\|B_2(u,v)\big\|_{L^2}^2\frac{dt}{t^{\alpha}}
&\lesssim \int_0^{r^2}\left\|\int_0^s e^{(s-h)\Delta}\Delta \Big((-\Delta)^{-\frac12}(I-e^{h\Delta})M(h,\cdot)\Big)dh\right\|_{L^2}^2\frac{dt}{t^\alpha}\\
&\lesssim\int_0^{r^2}\left\|\Big((-\Delta)^{-\frac12}(I-e^{s\Delta})M(s,\cdot)\Big)\right\|_{L^2}^2\frac{ds}{s^{\alpha}}.
\endalign
$$
Owing to $\sup_{s\in (0,\infty)}s^{-1}(1-\exp(-s^2))<\infty$, we conclude that $(-\Delta)^{-\frac12}(I-e^{s\Delta})$ is bounded on $L^2$ with operator norm $\lesssim \sqrt{s}$. This, plus the Cauchy-Schwarz inequality, gives
$$
\int_0^{r^2}\big\|B_2(u,v)\big\|_{L^2}^2\frac{dt}{t^{\alpha}}\lesssim r^{n-2\alpha}\|u\|_{(X_{\alpha;T})^n}^2\|v\|_{(X_{\alpha;T})^n}^2.
$$

Similarly for $B_3(u,v)$, we obtain

$$
\align
\int_0^{r^2}\big\|B_3(u,v)\big\|_{L^2}^2\frac{dt}{t^{\alpha}}&\lesssim\int_0^{r^2}\left\|(-\Delta)^\frac12 e^{t\Delta}\int_0^t M(s,\cdot)ds\right\|_{L^2}^2\frac{dt}{t^\alpha}\\
&\lesssim r^{4+n-2\alpha}\int_0^1\left\|(-\Delta)^\frac12 e^{\tau\Delta}\int_0^\tau |M(r^2\theta,r\cdot)|d\theta\right\|_{L^2}^2\frac{d\tau}{\tau^\alpha}.
\endalign
$$
Now, making a use of Lemma 3.2 we achieve
$$
\int_0^1\left\|(-\Delta)^\frac12 e^{\tau\Delta}\int_0^\tau |M(r^2\theta,r\cdot)|d\theta\right\|_{L^2}^2\frac{d\tau}{\tau^{\alpha}}\lesssim D(M;\alpha)\int_0^1\Big\|M(r^2\theta,r\cdot)\Big\|_{L^2}^2\frac{d\theta}{\theta^{\alpha}},
$$
where
$$
D(M;\alpha)=\sup_{\rho\in (0,1)}\rho^{-n}\int_0^{\rho^2}\int_{|w-x|<\rho} |M(r^2\theta,rw)|\tau^{-\alpha}dwd\tau\lesssim r^{-2}\|u\|_{(X_{\alpha;T})^n}\|v\|_{(X_{\alpha;T})^n}.
$$
Observe also that
$$
\int_0^1\|M(r^2\theta,r\cdot)\|_{L^2}^2\frac{d\theta}{\theta^{\alpha}}\lesssim r^{-2}\|u\|_{(X_{\alpha;T})^n}\|v\|_{(X_{\alpha;T})^n}.
$$
So, it follows that
$$
\int_0^{r^2}\|B_3(u,v)\|_{L^2}^2\frac{dt}{t^{\alpha}}\lesssim r^{n-2\alpha}\|u\|^2_{(X_{\alpha;T})^n}\|v\|^2_{(X_{\alpha;T})^n}.
$$
Adding the previous estimates on $B_j(u,v)$, $j=1,2,3$ together gives (3.5). 

Clearly, the boundedness of $B(\cdot,\cdot): (X_{\alpha;T})^n\times(X_{\alpha;T})^n\to (X_{\alpha;T})^n$ follows from (3.4) and (3.5). Furthermore, the case $T=\infty$ produces (i); and the other case $T\in (0,\infty)$ yields (ii). The proof is complete.
\enddemo

\vskip 0.6 true cm

\smallskip

\settabs\+xxxxxxxxxxxxxxxxxxxxxxxxxxxxxxxxx
&\cr\+ Department of Mathematics and Statistics&\ \ \ \cr\+ Memorial
University of Newfoundland &\ \ \
\cr\+St. John's, NL, A1C 5S7 &\ \ \ \cr\+ Canada &\ \ \
\cr\+Email: jxiao\@math.mun.ca&\ \ \
\cr\+&\cr


\Refs \widestnumber\key{WWW}

\ref \key{AdL}\by D. R. Adams and J. L. Lewis\paper On Morrey-Besov inequalities\jour {Studia\ Math.}\vol {74}\yr 1982\pages 169-182
\endref

\ref \key{AnC}\by M. Andersson and H. Carlsson\paper $Q_p$ spaces in strictly pseudoconvex domains\jour {J\ Anal.\ Math.}\vol {84}\yr 2001\pages 335-359
\endref


\ref \key{Ca}\by S. Campanato\paper Propriet\'a di inclusione per spazi di Morrey\jour {Ricerche Mat.}\vol {12}\yr 1963\pages 67-86
\endref

\ref \key{Ch}\by M. Christ\paper Lectures on Singular Integral Operators\jour{CBMS Regional Conference Series in Mathematics}\publ{Amer. Math. Soc., Providence, RI}\vol 77\yr 1990
\endref

\ref \key{CoMS}\by R. R. Coifman, Y. Meyer and E. M. Stein\paper Some new function spaces and their applications to harmonic analysis\jour {J.\ Funct.\ Anal.}\vol {62}\yr 1985\pages 304-335
\endref

\ref \key{CuY}\by L. Cui and Q. Yang\paper On the generalized Morrey spaces\jour {Siberian\ Math.\ J.}\vol {46}\yr 2005\pages 123-141
\endref

\ref \key{DX1}\by G. Dafni and J. Xiao\paper Some new tent spaces and duality theorems for fractional Carleson measures and $Q_\alpha(\rn)$\jour {J.\ Funct.\ Anal.}\vol {208}\yr 2004\pages 377-422
\endref

\ref \key{DX2}\by G. Dafni and J. Xiao\paper The dyadic structure and atomic decomposition of $Q$ spaces in several variables\jour{Tohoku\ Math.\ J.}\vol 57\yr 2005\pages 119-145
\endref

\ref \key{En}\by M. Englis\paper $Q_p$-spaces: generalizations to bounded symmetric domains\jour {preprint}\yr 2005
\endref

\ref \key{EsJPX}\by M. Ess\'en, S. Janson, L. Peng and J. Xiao\paper{$Q$ spaces of several real variables}\jour {Indiana\ Univ.\ Math.\ J.}\vol {49}\yr 2000\pages 575-615
\endref


\ref \key {FS}\by C. Fefferman and E. Stein\paper {$H^p$ spaces of several variables}\jour {Acta\ Math.}\vol {129}\yr 1972\page 137-193
\endref

\ref \key {JN}\by F. John and L. Nirenberg\paper {On functions of bounded mean oscillation}\jour {Comm.\ Pure\ Appl.\ Math.}\vol {18}\yr 1965\pages 415-426
\endref

\ref \key{Ka}\by T. Kato\paper {Strong $L^p$-solutions of the Navier-Stokes equation in $\Bbb R^m$, with applications to weak solutions}\jour {Math.\ Z.}\vol {187}\yr 1984\pages 471-480
\endref

\ref \key{KoTa}\by H. Koch and D. Tataru\paper {Well-posedness for the Navier-Stokes equations}\jour {Adv.\ Math.}\vol {157}\yr 2001\pages 22-35
\endref

\ref \key{La}\by V. Latvala\paper{On subclasses of $BMO(B)$ for solutions of quasilinear elliptic equations}\jour {Analysis}\vol 19\yr 1999\pages 103-116
\endref

\ref \key {Le}\by P. G. Lemari-Rieusset\paper {Recent Developments in the Navier-Stokes Problem}\jour {Boca Raton, Fla.: Chapman and Hall/CRC }\yr 2002
\endref

\ref \key {Lie}\by E. Lieb\paper {Sharp constants in the Hardy-Littlewood-Sobolev and related inequalities}\jour {Ann.\ Math.}\vol {118}\yr 1983\pages 349-374
\endref

\ref \key {LieL}\by E. Lieb and M. Loss\paper {Analysis}, { 2nd ed. Graduate Studies in Mathematics}\vol {14}\publ Amer. Math. Soc., Providence, Rhode Island\yr 2001.
\endref

\ref \key {LinS}\by H. Lindblad and C. D. Sogge\paper {On existence and scattering with minimal regularity for semilinear wave equations}\jour {J.\ Funct.\ Anal.}\vol {130}\yr 1995\pages 357-426
\endref

\ref \key{MS}\by V. Maz\'ya and T. Shaposhnikova\paper{On the Bourgain, Brezis, and Mironescu theorem concerning limiting embeddings of fractional Sobolev spaces}\jour {J.\ Funct.\ Anal.}\vol 195\yr 2002\pages 230-238
\endref


\ref \key{Pe}\by J. Peetre\paper{On the theory of $\Cal L_{p,\lambda}$}\jour {J.\ Funct.\ Anal.}\vol 4\yr 1969\pages 71-87
\endref

\ref \key{PeY}\by L. Peng and Q. Yang\paper{Predual spaces for Q spaces}\jour {Preprint}\yr 2005
\endref


\ref \key{PoS}\by S. Pott and M. Smith\paper{Paraproducts and Hankel operators of Schatten class via $p$-John-Nirenberg theorem}\jour {J. Funct. Anal.}\vol 217\yr 2004\pages 37-78
\endref

\ref \key {Ste}\by E. M. Stein\paper {Singular Integrals and Differentiability Properties of Functions}\publ {Princeton University Press}, Princeton, N.J.\yr 1970
\endref

\ref \key{Str1}\by R. S. Strichartz\paper{Bounded mean oscillation and Sobolev spaces}\jour {Indiana\ Univ.\ Math.\ J.}\vol 29\yr 1980\pages 539-558
\endref

\ref \key{Str2}\by R. S. Strichartz\paper{A Guide to Distribution Theory and Fourier Transforms}\publ{World Scientific}\yr 2003
\endref


\ref \key{Ta}\by M. E. Taylor\paper{Analysis on Morrey spaces and applications to Navier-Stokes and other evolution equations}\jour {Comm.\ Partial\ Differential\ Equations}\vol 17\yr 1992\pages 1407-1456
\endref

\ref \key{WuXi}\by Z. Wu and C. Xie\paper{$Q$ spaces and Morrey spaces}\jour{J.\ Funct.\ Anal.}\vol{201}\yr 2003\pages 282-297
\endref

\ref \key{X}\by J. Xiao\paper{A sharp Sobolev trace inequality for the fractional-order derivatives}\jour{Bull.\ Sci.\ Math.}\vol{130}\yr{2006}\pages 87-96 
\endref

\ref \key{Y}\by M. Yamazaki\paper{
The Navier-Stokes equation in various function spaces}\jour {Amer.\ Math.\ Soc.\ Transl.}\vol 204\yr 2001\pages 111-132
\endref
\endRefs
\enddocument